\documentclass[11pt,twoside]{amsart}

\usepackage{amssymb}

\usepackage[all]{xy}
\usepackage{epsfig}  
\usepackage{color}
\usepackage{pifont}

\title{Steinness of bundles with fiber a Reinhardt bounded domain}
\author{Karl Oeljeklaus \& Dan Zaffran} 
\thanks {karl.oeljeklaus@cmi.univ-mrs.fr , zaffran@math.cornell.edu\ }

\date{\today}

\newtheorem{dfn}{Definition}[section]
\newtheorem{thm}{Theorem}
\newtheorem{lem}[dfn]{Lemma}
\newtheorem{crl}[dfn]{Corollary}
\newtheorem{prop}[dfn]{Proposition}

\def\l{\ensuremath{\lambda}}
\def\linv{\ensuremath{\frac{1}{\lambda}} }
\def\pmlinv{\ensuremath{\frac{\pm 1}{\lambda}} }

\def\inv{^{-1}}

\def\e{\text{\it \Large e\! }}

\def\lone{\mbox{$\lambda _1$}}
\def\ltwo{\mbox{$\lambda_2$}}
\def\lthr{\mbox{$\lambda_3$}}

\def\Cvx{\ensuremath{log\ D}}

\def\A{\ensuremath{\mathcal A}}
\def\P{\ensuremath{\mathcal P}}
\def\S{\ensuremath{\mathcal S}}
\def\O{\ensuremath{\mathcal O}}

\def\forall {\text{for all }}

\def\N{\ensuremath{\mathbb N}}
\def\Q{\ensuremath{\mathbb Q}}
\def\R{\ensuremath{\mathbb R}}
\def\Z{\ensuremath{\mathbb Z}}

\def\C{\ensuremath{\mathbb C}}

\def\CP1{\ensuremath{\mathbb C \mathbb P^1}}
\def\Cstar{\ensuremath{\C^*}}

\def\Pn-1{\ensuremath{\P^{n-1}}}

\def\caseone{\ensuremath{{\mathbf[1]}}}
\def\casetwoC{\ensuremath{{\mathbf[2]_\C}}}
\def\casetwoR{\ensuremath{{\mathbf[2]_\R}}}
\def\casethree{\ensuremath{{\mathbf[3]}}}

\newcommand{\mysmallmatrix}[1]{\ensuremath{
                            \bigl[ \begin{smallmatrix}
#1 
\end{smallmatrix} \bigr]}}

\newcommand{\applic}[5]{
  \ensuremath{\begin{array}[t]{c@{}c@{}c@{}l}
           #1:\ & #2 & \rightarrow & #3 \\
               & #4 & \mapsto     & #5
         \end{array}
        }
}

\newcommand{\vect}[1]{ {\scriptstyle{ 
\left[
\begin{array}{c} 
#1
\end{array}
\right]               }}}

\begin{document}

\maketitle
\thispagestyle{empty}

\section{Introduction and notations}
\label{intro}

\noindent
Stein manifolds can be characterized by the fact that they holomorphically embed 
in $\C^N$ for some $N$. From that point of view it is obvious that 
if $F$ and $B$ are Stein manifolds, then the product $E=F\times B$ 
is Stein. More generally, take a fiber bundle $E$ with fiber $F$ and 
with basis $B$, which we shall denote by $E\xrightarrow{F}B$. 
Is such an $E$ necessarily Stein? That question was asked fifty 
years ago by J.-P. Serre (cf. \cite{Ser}), and is often referred to 
as ``the Serre problem'' in the literature. A counterexample 
with $F=\C^2$ was produced by H. Skoda in 1977 (cf. \cite{Skoda}). 
On the other hand, many interesting ``positive results'' have been 
obtained (see below). 

Following \cite{Pfl-Zwo}, we shall say that a Stein manifold $F$ is of 
class \S, or $F\in\S$ for short, when\footnote{Conveniently enough, that letter 
honors simultaneously Serre, Sibony, Siu, Skoda, Stehl\'e 
and Stein.}\,:
For any bundle $E\xrightarrow{F}B$ with $B$ Stein, $E$ is Stein. 

\smallskip\noindent 
K. Stein proved in \cite{Stein} that if $dim\ F=0$, then $F\in\S$.
Building on previous work by A. Hirschowitz, Y.T. Siu and N. Sibony 
(cf. \cite{Hirsch}, \cite{Sib} and \cite{Siu1}), N. Mok proved in \cite{Mok} 
that if $dim\ F=1$, then $F\in\S$. Skoda's above result can be stated as: 
$\C^2 \not\in \S$. 

\smallskip
In this paper we focus on the case where $F$ is a bounded 
domain $D\subset\C^n$ (``domain'' means connected open subset). 
Several results showing that large classes of 
bounded domains belong to $\S$ have been proved 
(cf. \cite{Siu2}, \cite{Steh} and \cite{Die-For}).
Nevertheless, in \cite{Coe-Loe}, G. C\oe ur\'e and J.-J. 
L\oe b produced, for each given 
$A\in SL_{2}(\Z)$ with $|trace\ A|>2$, 
a non Stein bundle $E_{CL}\xrightarrow{D_{CL}}\Cstar$ whose fiber 
$D_{CL}$ is a bounded Stein 
domain subset of $(\Cstar)^2$. Thus $D_{CL}\not\in \S$. 
Their $D_{CL}$ has the following properties (see picture p.\pageref{pictureDCL}):\\
-- $D_{CL}$ has the Reinhardt symmetry, i.e., it is 
invariant by the map $(z_1,z_2)\mapsto (\alpha_1 z_1,\alpha_2 z_2)$ 
for any complex numbers $\alpha_1$ and $\alpha_2$ of modulus $1$.\\
-- $D_{CL}$ has an automorphism $g$ of the form 
$g(z_1,z_2)=(z_1^{a}z_2^{b},z_1^{c}z_2^{d})$ with 
$\mysmallmatrix{a&b\\c&d}=A$. 

\smallskip
The second named author generalized their construction and gave a better 
understanding of those bundles. In \cite{Zaf}, a key point is the 
existence of a 
{\em $g$-equivariant} open dense embedding 
$D_{CL}\hookrightarrow \hat{D}_{CL}$, where 
$\hat{D}_{CL}\setminus D_{CL}$ is an 
infinite chain of rational curves. Roughly speaking, the non-Steinness 
of $E_{CL}$ is ``explained'' by what happens near that 
infinite chain.

\smallskip
Our goal here is to answer the following ``converse" question:\\  
Let $D\subset \C^n$, with $n=2$ or $3$, 
be {\em any} bounded Stein Reinhardt domain. 
Does $D$ belong to \S? In other words, does there exist a bundle 
$E\xrightarrow{D}B$ with $B$ Stein and $E$ non-Stein? The answer 
is contained in Theorem \ref{which_D_are_in_S} below. It reveals 
that in dimension two, C\oe ur\'e-L\oe b's examples $D_{CL}$ (for all 
$A\in SL_{2}(\Z)$) are essentially 
the only Reinhardt bounded domains not in $\S$: 
all other examples are provided 
by $g$-invariant subdomains of some $D_{CL}$. 
Moreover, it is easily checked that the interior of the 
closure in $\hat{D}_{CL}$ of such a 
subdomain still contains the infinite chain of curves, so 
from the point of view of \cite{Zaf}, it is natural 
that those domains do not belong to \S. Indeed, 
proofs in \cite{Coe-Loe} and \cite{Zaf} apply almost verbatim 
to show that they do not belong to \S. 
We shall see that both methods and results become 
more complicated in dimension three.
 
We also address the following question: given $D$ bounded and Reinhardt not in 
\S\, and B Stein, can we give a characterization of the 
Steinness of a bundle $E\xrightarrow{D}B$? For a 
two-dimensional $D$, we give in Theorem 
\ref{which_bundles_are_Stein} both a necessary criterion and a sufficient criterion. For a higher dimensional 
$D$, we give a partial answer in Theorem \nolinebreak 
\ref{which_bundles_are_Stein_dim_n}.

\medskip
{\em 
We work throughout the article in the complex category. In other words, 
all manifolds and maps we deal with are holomorphic. 
By the word  ``bundle" we mean a locally trivial holomorphic fiber bundle. 
We shall also use the notations $\O(E)$ for the set of holomorphic 
functions on $E$, and  $S^1$ and $\Delta$ will respectively denote the unit circle 
and the unit disk in \C.
 
NB: Most proofs are postponed until the end of the paper, in section 
\ref{proofs_of_lemmas}. 
}

\medskip
We shall make use of several known results about a given bundle 
$E\xrightarrow{D}B$ with $B$ Stein and $D\subset \C^n$ bounded and Stein. 
Namely:\\ 
\begin{itemize}
\item $E$ is a flat bundle (see \cite{Kaup} or \cite{Siu2}). That 
means $E$ can be defined by locally constant transition functions. 
All flat bundles can be constructed as follows:\\ 
Let $\rho$ be a morphism $\pi_1(B)\rightarrow Aut(D)$. 
Such a $\rho$ induces a $\pi_1(B)$-action on $D$. 
Denote by $\tilde{B}$ the universal cover of $B$, and consider 
the diagonal action of $\pi_1(B)$ on $\tilde{B}\times D$. 
Define 
$$E = \frac{\tilde{B}\times D}{\pi_1(B)}.$$ 
Then the projection map 
$\tilde{B}\times D\rightarrow \tilde{B}$ induces a map 
$E\rightarrow B$ that turns $E$ into a bundle with fiber $D$. 
\label{E_given_by_representation}
The {\em structural group} $G_{struct}(E)$ of $E$ is by 
definition the image of $\rho$. That definition is quite improper 
because $G_{struct}(E)$ does not only depend on the isomorphism 
class of $E$ as an $F$-bundle over $B$, but also on the $\rho$ 
chosen. That ``subtlety'' won't matter for our purposes... 
\item $E$ is holomorphically separable (see \cite{Siu2}). 
\item $E$ is Stein if $G_{struct}(E)$ is compact 
or contained in a compact group (see \cite{Sib} and \cite{Siu2}). 
\item $E$ is Stein if $G_{struct}(E)$ has  finitely many connected components 
(see \cite{Siu2}). 
\end{itemize} 
Given $g\in Aut(D)$, there is exactly one bundle $E\xrightarrow{D}\Cstar$ 
with monodromy $g$. Namely, $E=\frac{\C\times D}{\Z}$, where 
the \Z-action is the ``diagonal action" generated by 
$$\applic{\tilde{g}}{\C \times D}{\C \times D}
{\Big( z\ ; d\Big)}{\Big( z+1\ ; g(d)\Big).}
$$
We shall call that bundle the {\em complex suspension of $g$}. It has infinite cyclic 
$\A(E)$, generated by $g$. 

\medskip \noindent
{\em For simplicity, we introduce the following notations and results 
assuming $n=2$, but they all extend in the ``obvious" way to any $n\geq 2$.} 

\smallskip
By a well-known criterion for the Steinness of a Reinhardt domain 
(see \cite{Nara}), the map
$$\applic{\text{``$log$''}}{(\Cstar)^2}{\R^2}{(z_1,z_2)}{(log\ |z_1|, log\ |z_2|)}$$ 
induces a one-to-one correspondence between Stein Reinhardt domains 
of $(\Cstar)^2$ and open convex subsets of $\R^2$. 

Now take $D\subset(\Cstar)^2$ a {\em bounded} Stein Reinhardt domain. 
We shall denote by $\Cvx$ the image of $D$ by the above map. To make more 
explicit the one-to-one correspondence $$\Cvx \leftrightarrow\ D,$$
remark 
that $D$ can be recovered from \Cvx\ as the image of the ``tube'' 
$\Cvx + i \R^2 \subset \C^2$ by the map 
$\text{``$exp$''}:\ (w_1, w_2)\mapsto (z_1,z_2)=(\e^{w_1},\e^{w_2})$. 
Moreover \Cvx\ contains no affine line: otherwise 
$\Cvx + i \R^2$ would contain a copy of \C\ on which ``$exp$'' 
would restrict to a non-constant map from \C\ to $D$, contradicting 
the boundedness of $D$. 
By \cite{Zwo}, the converse statement 
holds: 
for a Stein Reinhardt $D$, if \Cvx\ contains no affine line, 
then $D$ is isomorphic to a bounded domain (we won't use that fact, 
though). 

\smallskip 
We denote by $Aut(D)$ the group of automorphisms of $D$. 
By \cite{Shim}, $$Aut(D)=Aut_{alg}^{\R}(D) \ltimes Aut_{0}(D),$$ 
where:\\
--- $Aut_{0}(D)$ is the connected component of the identity,\\ 
--- $Aut_{alg}(D)$ is the subgroup of $Aut(D)$ defined by: 
For each $g\in Aut_{alg}(D)$, there exist 
$A_g=\mysmallmatrix{a&b\\c&d} \in GL_2(\Z)$ 
and $\alpha_1,\alpha_2 \in \Cstar$ such that 
$$g(z_1,z_2)=(\alpha_1 z_1^a z_2^b, \alpha_2 z_1^c z_2^d).$$
(Given $g$, such $A_g$,  $\alpha_1$ and $\alpha_2$ are unique),\\
--- $Aut_{alg}^{\R}(D)$ denotes the subgroup of $Aut_{alg}(D)$ 
of all $g$'s with real and positive $\alpha_i$'s.   
Thus 
$Aut_{alg}(D) = Aut_{alg}^{\R}(D) \ltimes (S^1)^2$, and by lemma 
\ref{psi_iso},  $ Aut_{alg}^{\R}(D)$ is a discrete group.

\smallskip
For $g\in Aut_{alg}(D)$, we shall 
denote by $f_g$ the map 
$$\applic{f_g}{\Cvx}{\Cvx}{p}{A_g p + b_g,}$$ 
where $b_g=(log\ |\alpha_1|,\ log\ |\alpha_2|)$. 
The correspondence $g\leftrightarrow f_g$ is one-to-one 
between $Aut_{alg}^{\R}(D)$ and the group of affine 
automorphisms of \Cvx. 
Remark that $f_{g\inv}(p)=f_g\inv(p)=A_g\inv p - A_g\inv b_g$.
 
Define 
$$\A(D)=\{ A_g\ :\ g\in  Aut_{alg}^{\R}(D)\} \subset GL_2(\Z).$$
It is useful to think of  $\A(D)$ as ``the set of matrices that act on $D$".
 
For a given bundle $E\xrightarrow{D}B$, 
$$\A(E)=\Big\{ A_g\ :\ g\in  Aut_{alg}^{\R}(D) \cap 
\big(Aut_0(D).G_{struct}(E)\big)\Big\} \subset \A(D).$$
It is useful to think of  $\A(E)$ as ``the set of matrices that are used to build $E$".

For any group of matrices $\A$, we denote 
$$Spec_{\C}\ \A= \bigcup_{A \in \A} Spec_{\C}\ A.$$

\bigskip
We can now state the main results of this paper. They consist of 
the following theorems, and the geometric description (that follows 
from Theorem \ref{which_D_are_in_S}) given below. We point out 
the importance of the set $Spec_{\C}\ \A(D)$ to study whether 
a domain $D$ belongs to \S\ or not. 
Theorem \ref{which_D_are_in_S} below in the case $n=2$ 
is the main result of \cite{Pfl-Zwo}. Our proof for that case is simpler. 
\begin{thm} 
\label{which_D_are_in_S} 
A bounded Stein Reinhardt domain $D\subset \C^n$ with $n=2$ or $3$  
belongs to \S\ if and only if $Spec_{\C}\ \A(D) \subset S^1$.
\noindent
In other words: 
\begin{itemize}
\item For $n=2$, $D\not\in \S$ if and only if 
there exists $A\in \A(D)$ with 
$Spec_{\C}\ A=\{\l,\linv\}$, $\l \in \R\setminus \{+1,-1\}$.  
\item For $n=3$, 
$D\not\in \S$ if and only if (cf. Lemma \ref{cases_for_Spec}) 
there exists $A\in \A(D)$ such that either
\begin{itemize}
  \item[(a)] $Spec_\C\ A = \{\lambda_1, \lambda_2, \lambda_3 \}$ 
                     with $\lambda_i$'s pairwise distinct and real, or 
  \item[(b)] $Spec_\C\ A = \{1,\l,\linv\}$ with $\l \in \R\setminus \{+1,-1\}$. 
\end{itemize}
\end{itemize}
\end{thm}
\noindent 
In any manifold $B$, consider a free homotopy class $\gamma:S^1 \rightarrow B$. 
By a real-analytic approximation, followed by complexification, $\gamma$ can be realized 
as a holomorphic map $$\{ w\in \C : 1<|w|<m \} \rightarrow B$$ for some $m\in (1, +\infty]$. 
We define the {\em holomorphic width} of $\gamma$ to be the $\log$ of the supremum of all possible such $m$'s. 


\begin{thm} 
\label{which_bundles_are_Stein} 
Let $D\subset \C^2$ be a bounded Stein Reinhardt domain. 
Let $E$
be a bundle over a Stein basis $B$, and with fiber $D$. If $Spec_\C\ \A(E) \subset S^1$ then $E$ is Stein. 

Conversely, assume that $Spec_\C\ \A(E) \not\subset S^1$, i.e., $B$ contains a homotopy class $\gamma$ 
with associated monodromy $g \in Aut\ (D)$ such that 
$\rho>1$, where $\rho$ denotes the spectral radius of $A_g$. 
If the holomorphic width of $\gamma$ is strictly greater than $\frac{2\pi^2}{\log \rho}$, then $E$ is not Stein. 
\end{thm} 

\begin{thm} 
\label{which_bundles_are_Stein_dim_n} 
Let $D\subset (\Cstar)^n$ be a bounded Stein Reinhardt domain, with $n\geq 2$. 
Let $E$
be a bundle over a Stein basis, and with fiber $D$. If  $Spec_\C\ \A(E) \subset S^1$, 
then $E$ is Stein. 
\end{thm}

Theorem \ref{which_D_are_in_S} naturally leads to the following 
 
\smallskip \noindent
{\bf Question.}
{\em For an arbitrary $n$, does a bounded Stein Reinhardt domain 
$D\subset \C^n$ belong to \S\ if and only if $Spec_{\C}\ \A(D) \subset S^1$? }

\smallskip \noindent
We believe the answer to be ``yes''.

\medskip \noindent
{\bf Geometric description of domains not in \S}

\smallskip
{\em \bf  For n=2}. It follows from Theorem \ref{which_D_are_in_S} that if 
$D \not \in \S$, then $D$ has an automorphism $g$ of the form 
$g(z_1,z_2)\mapsto (z_1^{a}z_2^{b},z_1^{c}z_2^{d})$ with 
$A=\mysmallmatrix{a&b\\c&d}\in GL_{2}(\Z)$ and 
$Spec_\C\ \A(E) \not\subset S^1$. Up to replacing $A$ by $A^{2}$, 
we can assume that  $A\in SL_{2}(\Z)$ and 
$Spec_{\C}\ A=\{\l,\linv\}$, $\l \in (1, +\infty)$. 
Then no entry of $A$ is zero, therefore $g$ maps both axes to $(0,0)$, 
so $D$ must be a subset of $(\Cstar)^{2}$. Moreover the action of $g$ 
on $(\Cstar)^{2}$ is well-understood: 

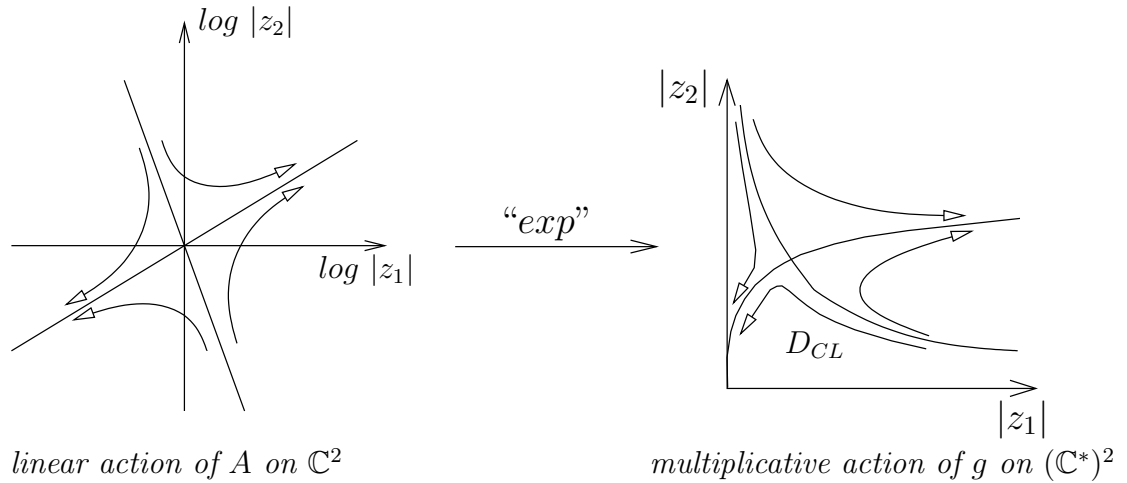
\begin{center} 
\hspace*{-7mm}

\label{pictureDCL}
\begin{picture}(0,0)%
\includegraphics{Quadrants.pstex}%
\end{picture}%
\setlength{\unitlength}{4144sp}%
\begingroup\makeatletter\ifx\SetFigFont\undefined%
\gdef\SetFigFont#1#2#3#4#5{%
  \reset@font\fontsize{#1}{#2pt}%
  \fontfamily{#3}\fontseries{#4}\fontshape{#5}%
  \selectfont}%
\fi\endgroup%
\begin{picture}(6166,2869)(754,-2369)
\put(4591,-2311){\makebox(0,0)[lb]{\smash{\SetFigFont{12}{14.4}{\familydefault}{\mddefault}{\updefault}{\color[rgb]{0,0,0}{\em multiplicative action of $g$ on $(\Cstar)^2$}}%
}}}
\put(766,-2311){\makebox(0,0)[lb]{\smash{\SetFigFont{12}{14.4}{\familydefault}{\mddefault}{\updefault}{\color[rgb]{0,0,0}{\em linear action of $A$ on $\C^2$}}%
}}}
\put(6661,-2041){\makebox(0,0)[lb]{\smash{\SetFigFont{14}{16.8}{\familydefault}{\mddefault}{\updefault}{\color[rgb]{0,0,0}$|z_1|$}%
}}}
\put(4636,-61){\makebox(0,0)[lb]{\smash{\SetFigFont{14}{16.8}{\familydefault}{\mddefault}{\updefault}{\color[rgb]{0,0,0}$|z_2|$}%
}}}
\put(2611,-1141){\makebox(0,0)[lb]{\smash{\SetFigFont{12}{14.4}{\familydefault}{\mddefault}{\updefault}{\color[rgb]{0,0,0}$log\ |z_1|$}%
}}}
\put(1891,344){\makebox(0,0)[lb]{\smash{\SetFigFont{12}{14.4}{\familydefault}{\mddefault}{\updefault}{\color[rgb]{0,0,0}$log\ |z_2|$}%
}}}
\put(3691,-871){\makebox(0,0)[lb]{\smash{\SetFigFont{14}{16.8}{\familydefault}{\mddefault}{\updefault}{\color[rgb]{0,0,0}$``exp"$}%
}}}
\put(5401,-1591){\makebox(0,0)[lb]{\smash{\SetFigFont{12}{14.4}{\familydefault}{\mddefault}{\updefault}{\color[rgb]{0,0,0}$D_{CL}$}%
}}}
\end{picture}

\label{Quadrants} 
\end{center} 

\noindent
The eigenspaces of $A$ split $\R^2$ into four quadrants, which correspond 
to four $g$-invariant ``quadrant domains" in $(\Cstar)^{2}$ that are depicted 
above. Any quadrant domain is sent isomorphically to a bounded domain by 
a well-chosen automorphism of $(\Cstar)^{2}$. The domain $D_{CL}$ appearing 
(as the fiber of a non Stein bundle) 
in \cite{Coe-Loe} is a quadrant domain of the matrix 
$\mysmallmatrix{2&1\\1&1}$.\\ 
Now, boundedness of $D$ implies that $D$ must be contained in 
a quadrant domain (otherwise $log$ $D$, being $A$-invariant, would contain an affine line). 

\smallskip
{\em Finally we conclude that for $n=2$, $D\not\in \S$ if and 
only if $D$ is a $g$-invariant Stein subdomain of a quadrant domain of 
some matrix $A\in SL_{2}(\Z)$. 

Such $D$'s are in one-to-one correspondence 
with convex $A$-invariant open subsets of a quadrant of $A$.} 

Basic examples are 
provided by the quadrants domains themselves, and also by taking $D$ such 
that $\log\ D$ is the interior of the convex hull of the $A$-orbit of a point 
in $\R^2$ not belonging to an eigenspace of $A$. 
 
\smallskip 
{\em \bf For n=3}. As opposed to the two-dimensional case, there are 
examples of domains not in \S, but not 
included in $(\Cstar)^3$: the product of the unit disk and any 
two-dimensional domain not in \S\ is easily checked not to be in \S.
Besides those trivial ``product examples'', domains 
invariant by a matrix that falls in case (b) of 
Theorem \ref{which_D_are_in_S} 
are some sort of a 
product with a two-dimensional domain $D'$ not in \S\ 
(see section \ref{product_case}).

Let $A$ be a matrix in $SL_3(\Z)$ satisfying condition (a) of 
Theorem \ref{which_D_are_in_S}, with 
eigenvectors $X_1, X_2$ and $X_3$. Examples of domains not in \S\ are: 
\begin{enumerate}
\item Any $D$ such that \Cvx\ is an octant delimited by the $X_i$'s. 
This example is the three-dimensional generalization of \cite{Coe-Loe} and \cite{Zaf}. 
\item Another basic example is given as follows: pick a point 
$p=x_1 X_1 + x_2 X_2 + x_3 X_3$ with $x_1x_2x_3\neq 0$. Then take 
$D$ such that \Cvx\ is the interior of the 
convex hull of $\{A^k p\ :\ k\in \Z \}$. 
\item In some cases, starting with a \Cvx\ from the above example, then 
reflecting it through a ``coordinate plane'', then taking the convex 
hull gives another convex $A$-invariant open subset of $\R^3$ that 
still corresponds to a {\em bounded} domain. 
{\em This produces examples not contained in 
an octant delimited by the $X_i$'s. These have no counterpart 
in the case of a two-dimensional $D$.} 
\end{enumerate} 

\smallskip
\noindent
{\bf Remark.} We've given a complete description of all domains not in \S\ 
for $n=2$, and a fairly precise description for $n=3$. But a complete 
answer should also tell when two domains in the list are actually 
isomorphic. 
For $n=2$, it actually follows from the above and from the literature 
about Inoue-Hirzebruch surfaces that considering 
``even Dloussky matrices'' (defined 
in \cite{Zaf}) already yields all examples of $D$'s not in \S.

\bigskip
\noindent
We shall make use of the following auxiliary results, proved 
in section \ref{proofs_of_lemmas}:

\begin{lem}
The only possible rational eigenvalues for 
a matrix $A\in GL_n(\Z)$ are $+1$ or $-1$. 
\label{only_plus_or_minus_one}
\end{lem}

\begin{lem}
Let $A\in GL_3(\Z)$. Then exactly one of the following holds:
\begin{itemize}
  \item[\caseone\: ]
                     $Spec_\C\ A \subset \{+1,-1\}$, 
  \item[\casetwoC]
                     $Spec_\C\ A = \{\pm 1,\l,\linv\}$ with 
                     $\l \in \C\setminus \{+1,-1\}$ and $|\l|=1$,
  \item[\casetwoR]
                     $Spec_\C\ A = \{\pm 1,\l,\pmlinv\}$ with 
                     $\l \in \R\setminus \{+1,-1\}$,  
  \item[\casethree\: ]
                     $Spec_\C\ A = \{\lambda_1, \lambda_2, \lambda_3 \}$ 
                     with $deg_\Q\ \lambda_i =3$ for $i=1\dots 3$. 
\end{itemize}
We shall write for short $GL_3(\Z)= \caseone \cup  \casetwoC 
\cup \casetwoR \cup \casethree$ (and that is a disjoint union). 
\label{cases_for_Spec}
\end{lem}

\begin{crl}
If $A\in GL_3(Z)$ doesn't have pairwise distinct eigenvalues, 
then $Spec_\C\ A \subset \{+1,-1 \}$.
\label{few_multiple_roots}
\end{crl}

\begin{lem}
The set of matrices $\A(D)$ is a subgroup of $GL_n(\Z)$, and the 
map $$\applic{\psi}{Aut_{alg}^{\R}(D)}{\A(D)}
                      {g}{A_g}$$
is a group isomorphism.
\label{psi_iso}
\end{lem}

\noindent
The following well-known result (see \cite{Hell} B(1)(b) p.47) 
is in some sense a generalization of 
B\'ezout's theorem (which corresponds to the case $n=2$). 
\begin{lem}
The vector $v=\vect{v_1\\ \vdots \\ v_n}\in \Z^n$ is unimodular 
(i.e., \mbox{$gcd(v_1,\dots,v_n)=1$}) if and only if there 
exists a matrix in $SL_n(\Z)$ of the form:
$$
\begin{bmatrix}
v_1     &    *    & \cdots &    *   \\
\vdots  & \vdots  &       & \vdots \\ 
v_n     &    *    & \cdots &    * 
\end{bmatrix}.
$$
\label{SLn_basis}
\end{lem}

\begin{lem}
\label{proj_Stein_Reinhardt}
Let $D\subset \C^n$ be a Stein Reinhardt domain, with $n\geq 2$. 
Denote $\pi(z_1, \ldots, z_n)=(z_2,\ldots, z_n)$. Then $\pi(D)$ 
is a Stein (Reinhardt) domain of $\C^{n-1}$. 
\end{lem}

\begin{lem}
Let $E\xrightarrow{D}B$ be a bundle with fiber a Reinhardt 
domain $D\subset (\Cstar)^n$. Assume there exists a non-zero $v\in\Z^n$ 
such that for all $A \in \A(E),\ Av=v$. Then there exists a ``quotient 
bundle'' $E'\xrightarrow{D'}B$ such that: $D'\subset (\Cstar)^{n-1}$;
$E'$ is Stein if and only if $E$ is Stein; and for any $A' \in \A(E')
\subset GL_{n-1}(\Z),\ A'$ is the automorphism of 
$\frac{\Z^n}{<v>_{\Z}}\approx \Z^{n-1}$ induced by some $A \in \A(E)$. 
In particular, 
$\ Spec_\C\ A(E') \subset Spec_\C\ A(E)$. 
\label{lowering_dimension}
\end{lem}

\noindent
The following result is due to Burnside (cf. \cite{Steinb} p.33).  

\begin{lem}
\label{Burnside}
Let $V$ be a finite-dimensional vector space over a field $F$. 
Let $S$ be a semi-group of endomorphisms of $V$ acting irreducibly 
on $V$. If the elements of $S$ have only finitely many different 
traces, then $S$ itself is finite. 
\end{lem}

\begin{lem}
Let $E\xrightarrow{D} B$ be a flat bundle and $H\subset G_{struct}(E)$ 
be a subgroup of finite index. Then there exists a flat bundle 
$E'\xrightarrow{D} B$ with $G_{struct}(E')=H$ and such that $E'$ 
is Stein if and only if $E$ is Stein. 
\label{finite_index_subgroup}
\end{lem}


\begin{lem}
Let $D\subset (\Cstar)^n$ be a bounded Stein Reinhardt domain, and 
$E\xrightarrow{D}B$ be a bundle with $B$ Stein. 
If $Spec_\C\ \A(E)$ is finite, then $E$ is Stein. 
\label{if_spec_finite}
\end{lem}

\begin{lem}
Let $\A$ be a subgroup of $GL_{n}(Z)$ with $n\geq 2$. 
If $Spec_\C\ \A\subset S^1$, then $Spec_\C\ \A$ is finite. 
\label{if_spec_subset_S1}
\end{lem}

\section{Proof of Theorem \ref{which_D_are_in_S} for $D\subset \C^2$.}
\label{dimension_two}
\noindent
The following two results prove Theorem \ref{which_D_are_in_S} for $n=2$. 

\begin{prop} \label{dim_two_in_Cstar}
Let $D\subset (\Cstar)^2$ be a bounded Stein Reinhardt domain. 
Then $D \in \S$ if and only if $Spec_\C\ \A(D) \subset S^1$.  
\end{prop}

\begin{proof}

\noindent {\em ``Only if" part: }
Suppose there exists $A\in \A(D)$ and $Spec_\C\ A \not\subset S^1$. 
Up to replacing $A$ by $A^2$, we can assume that $det\ A=1$. 
As $1\not\in Spec_\C\ A$, 
up to conjugation by an automorphism of $(\Cstar)^2$ , we can assume 
(see p.\pageref{Spec_does_not_contain_one} for more details) 
that $D$ has an automorphism of the form 
$g(z_1,z_2)=(z_1^a z_2^b,z_1^c z_2^d)$ with 
$A=\mysmallmatrix{a&b\\c&d}\in SL_2(\Z)$ and 
$Spec_\C\ A = \{ \l, \linv \}$, $\l \neq 1$. 
Let $E\xrightarrow{D}\Cstar$ be the complex suspension of $g$. 
Then 
the method of \cite{Coe-Loe} applies almost directly 
to show that $E$ is not Stein. The only (easy) extra work  
is to check that $log\ D$ contains an $\R$-orbit of 
the form $\{ A^t\mysmallmatrix{x\\y}\ :\ t\in\R \}$ for some 
$\mysmallmatrix{x\\y} \in log\ D$ (see section \ref{dim_three} for more details). 
The proof from \cite{Zaf} can be adapted 
as well. We conclude that $D\not\in\S$.\\ 

\medskip
\noindent {\em Conversely: }
Suppose that $Spec \A(D) \subset S^1$. 
Then by Lemma \ref{if_spec_subset_S1}, $Spec \A(D)$ is finite. 
Take a bundle $E\xrightarrow{D}B$ with $B$ Stein. 
As $Spec \A(E) \subset Spec \A(D)$, $E$ is Stein by lemma 
\ref{if_spec_finite}. This proves that $D\in \S$. 

\end{proof}

\begin{prop} \label{prop_dimension_two}\ \\
If $D\subset \C^2$ with $D\not\subset (\Cstar)^2$, then $D\in\S$ 
and $Spec_\C\ \A(D)\ \subset S^1$. 
\end{prop}
\begin{proof}
If $D\cap \{z_{i}=0\}$ is nonempty for $i=1,2$, then by 
\cite{Shim} p. 412, Lemma 1, $Aut_{alg}^{\R}(D)$ is finite.  
Therefore $Aut(D)$ has only finitely many connected components, so 
$D\in\S$ by \cite{Siu2}. As any $g\in Aut_{alg}^{\R}(D)$ has 
finite order, we also know that $Spec_{\C}\ \A(D) \subset S^1$. 

Else we can assume without loss of generality that 
$D\cap \{ z_{1}=0 \} \neq\varnothing$ and 
$D\cap \{ z_{2}=0 \}  = \varnothing$. Then any 
$g\in Aut_{alg}^{\R}(D)$ must have the form 
$g(z_{1},z_{2})=(\alpha_{1} z_{1} z_{2}^k, \alpha_{2} z_{2}^{\pm 1})$. 
As $D$ is bounded and stable by $g^n$ for any $n\in \Z$, we get $\alpha_{2}=1$, 
so the second component of $g^n(z_{1},z_{2})$ has modulus independent from $n$. 
Now, for the first component to remain bounded for any $(z_{1},z_{2})\in D$, 
we must have $\alpha_{1}=1$ and $k=0$. So $Aut_{alg}^{\R}(D)$ is finite 
and, as above, $D\in\S$ and $Spec_{\C}\ \A(D) \subset S^1$. 
\end{proof}

\section{Proof of Theorem \ref{which_D_are_in_S} for $D\subset \C^3$.}
\label{dim_three}

\subsection{For $D\subset (\Cstar)^3$, 
if $\A(D)\subset \caseone \cup \casetwoC$, then $D\in\S$.} 
\mbox{}\\
As $Spec_{\C}\ \A(D) \subset S^1$, it follows 
from lemmas \ref{if_spec_subset_S1} and \ref{if_spec_finite} 
that $D\in \S$. 


\subsection{For $D\subset (\Cstar)^3$, 
if $\A(D)\cap \casetwoR \neq \varnothing$, then $D\not\in\S$.} 
\label{product_case}
\mbox{}\\
Let $g\in Aut_{alg}(D)$ be such that $A=A_g \in \casetwoR$. 
Up to taking $A^2$ instead of $A$, we can assume that 
$Spec_\C\ A=\{1, \l, \linv \}$, with $\l \in\R\setminus\{+1, -1 \}$. 
As $1\in Spec_\Q\ A$, there exists $v\in\Q^3$ such that $Av=v$. 
Up to multiplying by some integer, we can assume that $v\in\Z^3$. 

Now define $E\xrightarrow{D}\Cstar$ as the complex suspension 
of $g$. 
Lemma \ref{lowering_dimension} applies here and yields 
a bundle $E'\xrightarrow{D'}\Cstar$ with fiber $D'\subset (\Cstar)^2$, 
and $\A(E')$ being generated by $A'\in SL_2(\Z)$ with 
$Spec_\C\ A'= \{\l, \linv \}$ and $\l\neq \pm 1$. 
We saw that $E'$ is not Stein in the proof of 
Proposition \ref{dim_two_in_Cstar}.
Therefore $D\not\in\S$. 

\subsection{For $D\subset (\Cstar)^3$, 
if $\A(D) \cap \casethree \neq \varnothing$, then $D\not\in\S$}
\label{generalized_CL}
\mbox{}\\
Let $g\in Aut_{alg}^{\R}(D)$ be such that $A=A_g \in \casethree$. 
We know that $g$ on $D$ corresponds on $\Cvx$ to 
$$f_g: \mysmallmatrix{x_1\\x_2\\x_3} 
\mapsto A\mysmallmatrix{x_1\\x_2\\x_3} + \mysmallmatrix{b_1\\b_2\\b_3}.$$ 
By Lemma \ref{cases_for_Spec}, $1\not\in Spec_\C\ A$ 
\label{Spec_does_not_contain_one}. Therefore 
$A-I$ is an invertible $3\times 3$ matrix, and it is easy to check that 
$p=-(A-I)\inv \mysmallmatrix{b_1\\b_2\\b_3}$ is a fixed point of $f_g$, and 
(therefore) the translation $T_p:\R^3\rightarrow \R^3$ of vector 
$p$ conjugates $f_g$ to $A$ (i.e., $T_p\inv f_g T_p = A$). 
It follows that up to conjugation by the $(\Cstar)^3$ 
automorphism $(z_1,z_2,z_3) \mapsto (\e^{p_1}z_1,\e^{p_2}z_2,\e^{p_3}z_3)$, 
we can assume that $g$ is such that 
$f_g \mysmallmatrix{x_1\\x_2\\x_3} = A\mysmallmatrix{x_1\\x_2\\x_3}$ 
for any $\mysmallmatrix{x_1\\x_2\\x_3}\in \R^3$. 

\begin{lem}
The eigenvalues $\l_1, \l_2, \l_3$ of $A$ are real. 
\end{lem}
\begin{proof}
As $\R^3$ is an odd-dimensional real vector space, we shall assume without 
loss of generality that $\l_1 \in \R$. We know from Lemma \ref{cases_for_Spec} 
that $\l_1 \neq 1$. Up to replacing $A$ by $A\inv$, we can assume that 
$\l_1 >1$.

Now we prove by contradiction that $\l_2$ is real. 
Suppose that $\l_2 \in \C\setminus\R$. As $det\ A=\pm 1$, we must have 
$|\l_2| <1$. Denote by $X_1$ a $\l_1$-eigenvector, and by $X_2$ and 
$X_3$ the imaginary and real parts of a (necessarily complex) 
$\l_2$-eigenvector. The matrix of $f_g$ with respect to the 
basis $\{X_1, X_2, X_3\}$ is:
$$
\begin{bmatrix}
\l_1 & 0 & 0 \\
0    & a & -b\\
0    & b & a 
\end{bmatrix},
$$
where $a+ib=\l_2$ (as $\l_2\not\in \R$, $b$ must be nonzero). 
As \Cvx\ is open, we can pick a point $q=x_1 X_1 + x_2 X_2 + x_3 X_3$ 
in \Cvx\ with $x_1\neq 0$. Then it is not hard to check 
that the convex hull of $\{A^k q\ :\ k\in \Z \}$ (which {\em is} in \Cvx) 
must contain the half-space $\{a_1 X_1 + a_2 X_2 + a_3 X_3\ :\ a_1 x_1>0 \}$. 
This is impossible because \Cvx\ contains no affine line.
So $\l_2$ and (therefore) $\l_3$ are real. 
\end{proof}

Up to taking $A^{\pm 2}$ instead of $A$, and renaming the eigenvalues, 
we can assume that $0<\lthr<\ltwo<1<\lone$.\\
Let $O^+$ denote the octant 
$\{ \sum_j a_j X_j/\ \forall j, a_j>0 \}$.
Again, we shall use coordinates with respect to the 
basis $\{X_1,X_2,X_3\}$, namely we write ``$(p_1,p_2,p_3)$'' 
instead of ``$\sum_j p_j X_j$''.  


The next two lemmas help us understand the shape of $\Cvx$: 

\begin{lem}
\label{shape_of_Cvx}
If $p=(p_1,p_2,p_3)\in \Cvx \cap O^+$ then, for any $t\geq 0,\ \\
(t+p_1,p_2,p_3)\in \Cvx \cap O^+$ and $(p_1,p_2,t+p_3)\in \Cvx \cap O^+$. 
\end{lem}
\begin{proof}
Recall that $0<\lthr<\ltwo<1<\lone$. Denote 
$(x_n,y_n,z_n):=A^np=(\lone^n p_1,\ltwo^n p_2,\lthr^n p_3)$.

When $n\rightarrow +\infty$, $y_n$ and $z_n$ tend to zero, 
whereas $x_n$ tends to $+\infty$. As a consequence, the 
angle between the line $(p,A^np)$ and the line $(p,p+(1,0,0))$ 
tends to zero. Now take a ball $B$ in $\Cvx$, centered at $p$. 
Fix $t\geq 0$. Because of the above facts, the convex hull of 
$\{A^np\} \cup B$ will contain $(t+p_1,p_2,p_3)$ as soon as $n$ 
is big enough. That proves the first part of the statement. 

When $n\rightarrow -\infty$, $x_n$ tends to zero, and 
$y_n$ and $z_n$ tend to $+\infty$, with $y_n/z_n$ tending to zero.
Therefore the angle between the line $(p,A^np)$ and the line 
$(p,p+(0,0,1))$ tends to zero. The second part of the statement 
follows as in the above case. 
\end{proof}

\begin{lem}
There exists $r\in \Cvx$ such that $\{A^tr, t\in \R \}$ is 
contained in $\Cvx$. 
\label{good_point_r}
\end{lem}
\begin{proof}
We shall make repeated use of Lemma \ref{shape_of_Cvx} and of the 
inequalities $\lone >1,\ \ltwo<1,\ \lthr <1$.
Take any $q=(q_1,q_2,q_3)\in \Cvx$.

\smallskip
{\em First case:} Assume that  $q \in O^+$. 
Then $V_1:=(\lone q_1,q_2,q_3)$ 
is in $\Cvx$, and therefore $V_2:=(\lone q_1,q_2,\lthr\inv q_3)$ 
is in $\Cvx$ as well. On the other hand, 
$Aq=(\lone q_1,\ltwo q_2,\lthr q_3)$ is in $\Cvx$, so 
$V_3:=(\lone q_1,\ltwo q_2, q_3)$ is in $\Cvx$, and so is 
$V_4:=(\lone q_1,\ltwo q_2, \lthr\inv q_3)$.

Multiplying the first components by $\lone$ we get four other 
points in $\Cvx$:\\  $V_5:=(\lone^2 q_1,q_2,q_3)$, 
$V_6:=(\lone^2 q_1,q_2,\lthr\inv q_3)$, 
$V_7:=(\lone^2 q_1,\ltwo q_2, q_3)$, and 
$V_8:=(\lone^2 q_1,\ltwo q_2, \lthr\inv q_3)$. 
%

Let $r:=V_2$. Remark that $V_7=Ar$. As $\Cvx$ is convex, it 
contains the parallelepiped $\P$ with vertices $V_1, V_2,\dots, V_8$. 
Thus, if $\lone q_1\leq x_1 \leq \lone^2 q_1$, 
$q_2\leq x_2 \leq \ltwo q_2$, and $\lthr\inv q_3\leq x_3 \leq q_3$, 
then $(x_1,x_2,x_3)$ belongs to $\P$, and therefore to $\Cvx$. 
It follows that for any $t$ between $0$ and $1$, 
$A^tr=(\lone^t \lone q_1,\ltwo^t q_2,\lthr^t \lthr\inv q_3)$ 
belongs to $\Cvx$. As $\Cvx$ is invariant by $A$, we get that $A^tr$ 
belongs to $\Cvx$ for any real $t$.

\smallskip
{\em Second case:} If  $q \not\in O^+$, take a linear automorphism 
$B$ of $\R^3$ such that for all $j$, $B(X_j)=\pm X_j$ 
(thus $A^tB=BA^t$ for any $t\in \R$), 
and such that $Bq \in O^+$. Now apply the first case to get in 
$B(\Cvx)$ an $\R$-orbit, which is sent to an $\R$-orbit in $\Cvx$ 
by $B\inv$. 
\end{proof}

\bigskip \noindent
We now construct a non-Stein bundle $E\xrightarrow{D}\Cstar$ 
by adapting the construction of G. C\oe ur\'e and J.-J. L\oe b. 
Lemma \ref{good_point_r} is an important ingredient for the 
construction of analytic disks as in \cite{Coe-Loe}.

We consider the  ``tube above $\Cvx$'':
$$ V:= \{ \sum_j u_j X_j\ /\ (Im\ u_1,Im\ u_2,Im\ u_3)\in \Cvx \}\  
\subset \C^3.$$
We consider the action of $G_\R =\R \ltimes \R^3$ on 
$\C \times \C^3$ given by 
 $$ \begin{array}[t]{c@{}c@{}c}
  \C \times \C^3\ \times\ G_\R & \rightarrow & \C \times \C^3 \\
      \big( (z, \sum_j u_jX_j )\ ,\ (t,X) \big)  
    & \mapsto     & \big( t+z, X+ \sum_j \lambda_j^t u_j X_j \big).
         \end{array}$$

Define $\Omega:=\C \times V$. By construction, $\Omega$ is 
invariant by the action of $G_\Z := \Z \ltimes \Z^3 \subset G_\R$. 

Now take $E:=\Omega / G_{\Z}$ (cf. \cite{Coe-Loe}). Remark that $E$ 
is just the complex suspension of $g$.
We shall prove 
that $E$ is not Stein by constructing a family of analytic disks 
$\{ d_{R}: \bar{\Delta} \rightarrow E\}_{R>1}$ such that 
$d_{R}(0)$ tends to infinity in $E$ (topologically) when 
$R\xrightarrow{>}1$, whereas 
$d_{R}(\partial \bar{\Delta})$ remains in a compact subset of $E$ 
independent from $R$. Here $\Delta$ (resp. $\bar{\Delta}$) denotes 
the open (resp. closed) unit disk in \C, and $\partial \bar{\Delta}$  
denotes their boundary. 

Fix a real number $R>1$. In what follows, it is useful to keep in mind that the maps $f$, 
$g_1, g_2$, $g_3$ and 
$\tilde{d}$ depend on $R$.\\ 
Define $f:\Delta \rightarrow \C$ by 
$$f(z)=log\ i\frac{R+z}{R-z}.$$
Then  $0<Im\ f<\pi$, and $f$ is continuous on $\bar{\Delta}$.
 
Thus we can define $g_1, g_2$ and $g_3$ as holomorphic 
functions $\Delta \rightarrow \C$ such that on $\partial \bar{\Delta}$, 
$Im\ g_j= \lambda_j^{\mu_1\inv Re\ f}$ for $j=1,2,3$, where 
$\mu_1$ denotes $log\ \lone$. 

Define a map $\tilde{d}: \Delta \rightarrow \C\times\C^3$ by 
$$\tilde{d}(t):= (\mu_1\inv f(t), \sum_j g_j(t) r_j X_j),$$ 
where $r=(r_1,r_2,r_3)$ is given by Lemma \ref{good_point_r}. Notice 
that $\tilde{d}$ is continuous on $\bar{\Delta}$. 
\begin{lem}
For any $t\in\bar{\Delta}, \ \tilde{d}(t) \in \Omega$. 
\end{lem}
\begin{proof}
Recall that $\tilde{d}(t) \in \Omega$ if and only if 
$(r_1Im\ g_1(t),r_2Im\ g_2(t),r_3Im\ g_3(t))\in \Cvx$. 

First assume that $t\in \partial \bar{\Delta}$.\\  
Denote $(Im\ g_1 (t) r_1,Im\ g_2(t) r_2, Im\ g_3(t) r_3)$ by 
$\tilde{d}'(t)$. By definition of the $g_j$'s, 
$\tilde{d}'(t)= 
(\lambda_1^{\alpha} r_1,\lambda_2^{\alpha} r_2, \lambda_3^{\alpha} r_3)$
for some $\alpha \in \R$, thus $\tilde{d}'(t)$ is on the $\R$-orbit of $r$, 
which is in $\Cvx$ 
by the Lemma \ref{good_point_r}. 

Now for $s\in \Delta$, the componentwise harmonicity of 
$\tilde{d}'$ allows us to write in vector form 
$$ 
\tilde{d}'(s)=\int_{S^1} \tilde{d}'(t) P_s(t),
$$
where $P_s$ is the Poisson kernel at $s$. 
As $\int_{S^1} P_s =1$, the above shows that 
$\tilde{d}'(s)$ is in the convex hull of $\tilde{d}'(S^1)$, 
therefore is in $\Cvx$.
\end{proof}

\medskip
We define an analytic disk $d_R:\Delta \rightarrow E$ by 
$$d_R=q_{\Z}\circ\tilde{d},$$ 
where $q_{\Z}$ is the quotient map $\Omega \rightarrow E$. 
Remark that $d$ is continuous on $\bar{\Delta}$. 

\begin{lem}
There exists a compact set $K$ in $E$, independent from $R$, 
such that $d_R(\partial\bar{\Delta}) \subset K$. 
\end{lem}
\begin{proof}
It is not hard to check that $\tilde{d}(t)$ is equal up to action by $G_{\R}$ to
$$\big(i \mu_1\inv Im\ f(t),i \sum_j r_j \lambda_j^{-\mu_1\inv Re\ f(t)}
Im\ g_j(t) X_j\big).$$
When $t$ is in $\partial\bar{\Delta}$, the above becomes 
$$\big(i \mu_1\inv Im\ f(t),i \sum_j r_j X_j\big).$$ 

As $0<Im\ f< \pi$, there exists a compact set $\tilde{K}\subset\Omega$ 
independent from $R$, such that 
$\tilde{d}(\partial\bar{\Delta}) \subset \tilde{K}.G_{\R}$. 
Then $d_R(\partial\bar{\Delta})$ is included in 
$K=q_{\Z}(\tilde{K}).(G_{\R}/G_{\Z})$, which is compact because 
$G_{\R}/G_{\Z}$ is compact.
\end{proof}

On the other hand, 
$$\tilde{d}(0)= (\mu_1\frac{\pi}{2}, \sum_j g_j(0) r_j X_j),$$
and $Im\ g_1(0)
=\int_{S^1} Im\ g_1
= \int_{S^1} \lambda_1^{\mu_1\inv Re\ f}
=\int_{S^1} \e^{Re\ f}$, which tends to $+\infty$ when 
$R$ tends to $1$ (cf. \cite{Coe-Loe}).\\ 
By construction of $E$, a point $q_{\Z}\big( z\ ; (z_1, z_2, z_3)\big)$ 
in $E$ tends to infinity (topologically) if and only if either 
$Im\ z$ or $Im\ (z_1, z_2, z_3)$ tends to infinity. Therefore 
when $R$ tends to $1$, 
$d_R(0)=q_{\Z}(\tilde{d}(0))$ tends to infinity in $E$. 
By the maximum principle (for functions $\Delta \rightarrow \C$), 
any function on $E$ will be bounded on the sequence $\{ d_{1+1/n}(0)\}_{n\in\N}$, 
which proves that $E$ is not holomorphically convex, therefore $E$ is not Stein.

Conclusion: $D\not\in\S$.

\subsection{The case of $D\not\subset(\Cstar)^3$.} 
\label{case_C_three_not_Cstar_three}
\mbox{}\\
Up to action on $\C^3$ by an algebraic automorphism, 
we can assume that $D$ is in Shimizu's 
normal form (in his notations, that means $D=\tilde{D}$; see \cite{Shim} p. 410). 
Remark: That realization of $D$ may not be bounded.  

\smallskip
\begin{prop}
If $D\subset\C^3$ intersects several coordinate hyperplanes, then 
$Spec_{\C}\ \A(D) \subset S^1$ and $D\in\S$. 
\end{prop}
\begin{proof}\ 

{\em First step.}
If $D\cap {z_{i}=0}$ is nonempty for $i=1,2$ and $3$, then by 
\cite{Shim} p. 412, Lemma 1, $Aut_{alg}^{\R}(D)$ is finite.  
Therefore $Aut(D)$ has only finitely many connected components, and 
$D\in\S$ by \cite{Siu2}. As any $g\in Aut_{alg}^{\R}(D)$ has 
finite order, we also know that $Spec_{\C}\ \A(D) \subset S^1$. 

\smallskip {\em Second step.}
Without loss of generality, we can now assume that 
$D\cap \{z_{1}=0\}$ and $D\cap \{z_{2}=0\}$ are both nonempty and 
$D\cap \{z_{3}=0\}$ is empty. (In Shimizu's notations, we must have 
$\sum_{i=1}^s n_{i} =2,\ t=s+1,\ n_{t}=1$.)\\
Take $g\in Aut_{alg}^{\R}(D)$. Then $g^2$ must 
be of the form: $$g^2(z_{1},z_{2},z_{3})= 
(\alpha_{1} z_{1} z_{3}^{\mu_{1}},
\alpha_{2}z_{2} z_{3}^{\mu_{2}},
\alpha_{3}z_{3}).$$ 
In particular, $Spec_{\C}\ \A(D) \subset S^1$.  
As $proj_{3}(D)$ is bounded (in Shimizu's notations, projection 
of $D$ on $\C^{n_{s+1}} \times \cdots \times \C^{n_t}$ is always bounded) 
and $g^2$-invariant, $\alpha_{3}=1$. Now, for any 
$(z_{1},z_{2},z_{3}) \in D,\ proj_{3}\inv(z_{3})$ is bounded and $g^2$-invariant, 
so $\alpha_{1}=\alpha_{2}=1$ and $\mu_{1}=\mu_{2}=0$. It follows that 
$Aut_{alg}^{\R}(D)$ is finite. Therefore, again, $D\in\S$. 
\end{proof}

\medskip \noindent
{\em We can now assume (until the end of section \ref{case_C_three_not_Cstar_three}) 
that $D\cap \{z_1=0\}$ is nonempty, 
and $D\cap \{z_{2}=0\}$ and $D\cap \{z_{3}=0\}$ are both empty. }

\medskip\noindent
Then any $g\in Aut_{alg}^{\R}(D)$ must 
be of the form: $$(*)_{alg}\quad g(z_{1},z_{2},z_{3})= 
(\alpha_{1} z_{1} z_{2}^{\mu_{2}} z_{3}^{\mu_{3}},
\alpha_{2} z_{2}^{a} z_{3}^{b},
\alpha_{3} z_{2}^{c} z_{3}^{d}),$$ 
with $A'_{g}=\mysmallmatrix {a&b\\c&d} \in GL_2(\Z)$. 
Denote by $\A '$ the group 
$\{A'_{g} \subset GL_2(\Z) \ :\  g\in Aut_{alg}^{\R}(D) \}$. 
Remark that $Spec_{\C}\ \A(D) = \{ 1 \} \cup Spec_{\C}\ \A '$, so 
$Spec_{\C}\ \A(D) \subset S^1$ if and only if  $Spec_{\C}\ \A ' \subset S^1$.

\begin{prop}
If $Spec_{\C}\ \A(D) \not\subset S^1$, then $D\not\in\S$. 
\end{prop}
\begin{proof}

Take $g\in Aut_{alg}^{\R}(D)$ with $Spec_{\C}\ A_g \not\subset S^1$. 
By $(*)_{alg}$, $g$ induces an automorphism $g'$ of $D'= D\cap \{ z_1=0 \}$, 
and $Spec_{\C}\ A'_{g'}=\{ \l, \linv \}$ with $\l \in\R\setminus\{+1, -1 \}$. 

Define bundles $E\xrightarrow{D}\Cstar$ and $E'\xrightarrow{D'}\Cstar$ 
as complex suspensions of $g$ and $g'$ respectively. The map $D' \hookrightarrow D$ 
induces a map $E' \hookrightarrow E$ realizing $E'$ as a closed submanifold. 
By the results of section \ref{dimension_two}, 
we know that $E'$ 
is not Stein, so $E$ can't be Stein. 

Conclusion: $D\not\in\S$. 
\end{proof}

\medskip \noindent
{\em Therefore we shall assume until the end of section \ref{case_C_three_not_Cstar_three} 
that }
$$ Spec_{\C}\ \A(D) \subset S^1.$$

\smallskip
Our assumptions ($z_1$ vanishes somewhere on $D$ but $z_2$ and $z_3$ don't) 
imply by \cite{Shim} that there can be three cases: 
(Shimizu's $n_{s+1}, \ldots, n_{t}$ must all be $1$, otherwise 
by \cite{Shim} p. 411 we would have an $U(2)$-action on $D$, 
and $D$ would intersect at least two coordinate hyperplanes.) 
\begin{enumerate}
\item \label{Os}
$D$ is bounded, and $Aut_{0}(D)$ is only $(S^1)^3$. 
(This is Shimizu's $D=\tilde{D}_2$ case.)
\item \label{Qs}
$D=\Big\{  (z_{1},z_{2},z_{3})\ :\      z_1\in\C, 
(z_2 \e^{q_2 |z_1|^2},  z_3 \e^{q_3 |z_1|^2}) \in D\cap \{z_1=0 \}  \Big\}$, 
with $q_{2}, q_{3}$ non-negative real numbers and $q_{2}q_{3}\neq 0$. 
\item \label{Ps}
$D=\Big\{  (z_{1},z_{2},z_{3})\ :\      |z_1|<1, 
(z_2 (1- |z_1|^2)^{-q_2/2},  z_3 (1- |z_1|^2)^{-q_3/2}) \in D\cap \{z_1=0 \}  \Big\}$, 
with $q_{2}, q_{3}$ non-negative real numbers (Shimizu's notations for 
$q_{2}$ and $q_{3}$ are $p_1^{2}$ and $p_1^{3}$). 
\end{enumerate}

{\bf If $D$ is in case (\ref{Os})}\\
Then any $g$ in $Aut\ (D)$ has the form $(*)_{alg}$ with 
$\alpha_1, \alpha_2, \alpha_3$ in $\Cstar$.\\
Denote $P_{23}(z_1,z_2,z_3)=(z_2,z_3)$, $D_{23}= P_{23}(D)$, and 
$\hat{D}= \{ (z_1,z_2,z_3)\ :\ z_1 \in \C, (z_2,z_3) \in D_{23} \}$.\\ 
Take $E\xrightarrow{D}B$ with Stein $B$.
The $Aut(D)$-action on $D$ induces actions on $\hat{D}$ and $D_{23}$ such 
that the 
sequence of maps $D\hookrightarrow \hat{D} \twoheadrightarrow D_{23}$ 
is $Aut\ (D)$-equivariant. Then we get a sequence of maps of bundles 
above $B$:  
$$\xymatrix{
E\ \ar[dr]_{D} \ar@{^{(}->}[r]  & \hat{E} \ar@{->>}[r] \ar[d]_{\hat{D}} & 
     E_{23} \ar[dl]^{D_{23}}\\
 & B &               
}.$$

Moreover the map $\hat{D} \twoheadrightarrow D_{23}$ turns 
$\hat{D}$ into a (trivial) line bundle over $D_{23}$, 
and that bundle structure is preserved by the $Aut\ (D)$-action 
on $\hat{D}$. 
Therefore $\hat{E}$ is a line bundle above $E_{23}$.\\ 

Now, as $Spec_{\C}\ \A(E_{23}) \subset Spec_{\C}\ \A(D) \subset S^1$, 
it follows from section \ref{dimension_two} that $E_{23}$ Stein. 
Therefore  $\hat{E}$ is Stein by \cite{Mok}. By local triviality 
of our bundles, $E$ is a locally Stein open subset of $\hat{E}$, 
so $E$ is Stein by Docquier-Grauert theorem.
 
Conclusion: $D\in \S$.

\medskip
{\bf If $D$ is in case (\ref{Qs}) or (\ref{Ps})}\\
The huge formulas in  \cite{Shim} p. 411 tell us 
that in cases (\ref{Qs}) and (\ref{Ps}), any $g_{0}\in Aut_{0}(D)$ is of the form: 
$$(*)_{0}\quad g_{0}(z_{1},z_{2},z_{3})= 
\Big( h_1(z_1) ,
z_{2} h_2(z_1),
z_{3} h_3(z_1) \Big),$$
where $h_1$ can be any automorphism of $\C$ in case (\ref{Qs}) and 
any automorphism of the unit disk in case (\ref{Ps}), and $h_2$ and $h_3$ \
are two holomorphic functions determined by $h_1$. 

\begin{lem}
For any $g\in Aut_{alg}^{\R}(D)$, the integers $\mu_1$ and $\mu_2$ 
in $(*)_{alg}$ are both zero. 
\end{lem}
\begin{proof}
Take $g\in Aut_{alg}^{\R}(D)$ (written as in $(*)_{alg}$), and 
$g_0 \in Aut_0(D)$ (written as in $(*)_0$) with $h_1$ not a 
homogeneous function.\\ 
Direct computation yields:
\def\azX{\alpha_1\inv z_1 X}

\medskip\noindent
$gg_0g\inv(z_1,z_2,z_3)=
\Big(\alpha_1 h_1(\azX) X\inv h_2^{\mu_2}(\azX)h_3^{\mu_3}(\azX)\ ,$

\hfill $  z_2 h_2^{a+c}(\azX)\ ,\ 
 z_3 h_3^{b+d}(\azX)
\Big)$

\medskip\noindent
with $X= \alpha_2^{\mu_2 d  -\mu_3 c} 
         \alpha_3^{-\mu_2 b +\mu_3 a}
         z_2^{-\mu_2 d +\mu_3 c}
         z_3^{\mu_2 b  -\mu_3 a}$. As $Aut_0(D)$ is a normal 
subgroup 
of $Aut(D)$, $gg_0g\inv$ belongs to $Aut_0(D)$, so the above 
expression must be of the form $(*)_0$. In particular, 
the first component of $gg_0g\inv(z_1,z_2,z_3)$ does not depend 
on $z_2$ or $z_3$. As $h_1$ is not homogeneous, $X$ itself 
must be independent of $z_2$ or $z_3$. Therefore $\mu_2=\mu_3=0$.
\end{proof} 

\noindent
Now $(*)_{alg}$ can be stated more precisely as:
$$ (*)_{alg}' \quad g(z_1,z_2,z_3)=
(\alpha_{1} z_{1},
\alpha_{2} z_{2}^{a} z_{3}^{b},
\alpha_{3} z_{2}^{c} z_{3}^{d}).
$$

\noindent
As $Spec_{\C}\ \A(D) \subset S^1$, $Spec_{\C}\ \A ' \subset S^1$. 
Therefore, by lemmas  \ref{if_spec_subset_S1} and \ref{Burnside}, there are two possibilities: 

\smallskip
1. $\A '$ is finite.\\  
Then by $(*)_{alg}'$, $\A(D)$ is finite, so $D\in \S$ by \cite{Siu2}. 

\smallskip 
2. There exists $v=(p,q) \in\Q^2$ such that 
for any $A' \in \A ', A'v=v$.\\ 
Up to division or multiplication by some integer, we 
can assume that $v$ belongs to $\Z^2$ and is unimodular. By Lemma \ref{SLn_basis},  
there exists $(p',q')$ 
such that $\mysmallmatrix{p&p' \\ q&q'}$ 
belongs to $SL_{2}(\Z)$. Consider the automorphism 
$$\applic{h}{\C\times (\Cstar)^2}{\C\times (\Cstar)^2}
                   {(z_{1},z_{2},z_{3})}
                   {(z_{1},z_{2}^p z_{3}^{p'},z_{2}^q z_{3}^{q'}).}$$
We replace $D$ with $h(D)$. The automorphisms get conjugated by $h$, so
any $g_{0}$ and $g$ in $Aut_{0}(D)$ and $Aut_{alg}^{\R}(D)$ now take the forms: 

$$(*)_0' \quad g_0(z_{1},z_{2},z_{3})=
(h_1(z_1),z_2 h_2^{q'} h_3^{-p}(z_1),z_3 h_2^{-q}h_3^{p}(z_1)),$$
$$(*)_{alg}'' \qquad g(z_{1},z_{2},z_{3})=
(\alpha_1 z_1, \alpha_{2} z_{2}, \alpha_3 z_3 z_2^{\star}).$$

Take a bundle $E\xrightarrow{D}B$ with $B$ Stein.\\ 
Denote $P_{12}(z_1,z_2,z_3)=(z_1,z_2)$, $P_1(z_1,z_2)=z_1$, 
$D_{12}= P_{12}(D)$,\\
 $D_1= P_1(D_{12})$, 
$\hat{D}=\{ (z_1,z_2,z_3)\ :\ (z_1,z_2)\in D_{12}, z_3 \in \C \}$ and\\ 
$\hat{D}_{12}=\{ (z_1,z_2)\ :\ z_1\in D_1, z_2 \in \C \}$. 
These domains are Stein by Lemma \ref{proj_Stein_Reinhardt}.\\ 
The forms  $(*)_0'$ and $(*)_{alg}''$ show that we get an 
induced $Aut(D)$-action on $D_{12}$, $D_1$, $\hat{D}$ and $\hat{D}_{12}$ 
such that the sequence of maps 
$D\hookrightarrow \hat{D} \twoheadrightarrow D_{12} 
\hookrightarrow \hat{D}_{12} \twoheadrightarrow D_1$ 
is equivariant. So we get a corresponding sequence of maps of 
bundles above $B$:
$$\xymatrix{
E\ \ar[drr]_{D} \ar@{^{(}->}[r]  & 
\hat{E} \ar@{->>}[r] \ar[dr] & 
E_{12} \ar[d] \ar@{^{(}->}[r]  &
\hat{E}_{12} \ar@{->>}[r] \ar[dl] &
E_1 \ar[dll]^{D_1} \\
 & & B & &               
}.$$
In analogy with case (\ref{Os}) above, we see that 
$\hat{E}$ and $\hat{E}_{12}$ are line bundles above $E_{12}$ 
and $E_1$ respectively. Now, as $D_1$ is one dimensional, 
$E_1$ is Stein by \cite{Mok}. The same reasoning as in 
case (\ref{Os}) shows that $E$ is Stein.

Conclusion: $D\in\S$.


\section{Postponed proofs}
\label{proofs_of_lemmas}

{\bf Theorem \ref{which_bundles_are_Stein} p.\pageref{which_bundles_are_Stein}:}
\begin{proof}
During the proof Theorem \ref{which_D_are_in_S}, we already proved 
that if $Spec_{\C}\ \A(E)\subset S^1$, then $E$ is Stein for $D$ of dimension $2$ or $3$. 

\smallskip\noindent 
Conversely, 
assume that there exist $\varepsilon_1 >0$, 
an annulus $U\subset \C$ of modulus $m = \frac{2\pi^2}{\log \rho} + \varepsilon_1$ and
$\varphi: U \rightarrow B$
with associated monodromy $g \in Aut\ (D)$, such that 
$\rho>1$, where $\rho$ denotes the spectral radius of $A_g$. 
As $Spec_{\C}\ A_g \not\in S^1$, it follows by Prop. \ref{prop_dimension_two} that 
$D\subset (\Cstar)^2$ and that $D$ is an $A_g$-invariant subdomain 
of a ``quadrant domain" (or ``C\oe ur\'e-L\oe b's domain"; cf. introduction).

Then the pulled-back bundle 
$E'=\varphi^*(E)$ has monodromy $g$, 
therefore $E'$ is just the restriction to the annulus $U$ of the complex suspension 
bundle of $g$ (i.e., C\oe ur\'e-L\oe b's bundle built from the matrix $A_g$).

The universal cover of $U$ can be realized as 
$$\{w\in\C\ :\ -\varepsilon_2 <Im\ w< \frac{\pi}{\log \rho} + \varepsilon_2 \}$$ 
with deck automorphism $w\mapsto w+1$ and $\varepsilon_2>0$. 
As noticed in Prop. \ref{dim_two_in_Cstar}, the method of \cite{Coe-Loe} can be directly 
adapted to the bundle $E'$. One gets a sequence of holomorphic disks 
$d'_n:\Delta \rightarrow E'$ whose boundaries are contained in a compact $K$ independent of $n$, 
and such that $d'_n(0)\rightarrow \infty$ topologically in $E'$. Moreover, there is 
a fiber of $E'$ that contains $d'_n(0)$ for all $n$. 



Take a bundle map $\tilde{\varphi}: E' \hookrightarrow E$ over $\varphi$ that restricts 
to an isomorphism on each fiber. 
The holomorphic disks in $E$ defined by $d_n= \tilde{\varphi} \circ d'_n$ 
have boundaries in the compact set $\tilde{\varphi}(K)$, and centers that tend to infinity.
This shows that $E$ is not Stein. 
\end{proof}

{\bf Theorem \ref{which_bundles_are_Stein_dim_n} 
                      p.\pageref{which_bundles_are_Stein_dim_n}:}
\begin{proof}
The theorem follows from lemmas \ref{if_spec_subset_S1} and \ref{if_spec_finite}. 
\end{proof}

{\bf Lemma \ref{only_plus_or_minus_one} p.\pageref{only_plus_or_minus_one}:}
\begin{proof}
Let $\l \in \Q\cap Spec_{\C}\ A$. Denote by $P_A$ the characteristic 
polynomial of $A$. Then $P_A(\l)=0$, but as $P_A$ is monic and 
has integer coefficients, any rational root of $P_A$ must be in $\Z$. 
Thus $\l \in \Z$. As $\linv\in Spec_\C\ A\inv$, we obtain 
similarly $\linv \in\Z$. Therefore $\l = +1$ or $\l=-1$.
\end{proof}

{\bf Lemma \ref{cases_for_Spec} p.\pageref{cases_for_Spec}:}
\begin{proof}
Denote $Spec_\C\ A = \{\lambda_1, \lambda_2, \lambda_3\}$. 
Let $P_A$ be the characteristic polynomial of $A$. Let 
$P_{\l_1}\in \Q[X]$ be the minimal polynomial of $\lambda_1$ over $\Q$. 

\begin{itemize}
\item If $deg\ P_{\l_1}=1$. Then $\lambda_1\in\Q$, so by lemma 
\ref{only_plus_or_minus_one}, $\l_1=\pm1$. Therefore $P_A$ 
can be written $P_A\overset{(*)}{=}(X-\lambda_1)Q$ with $Q\in \Q[X]$ 
and $deg\ Q=2$.

\begin{itemize}
\item If $\lambda_2=\pm 1$, $det\ A=1$ implies $\lambda_3=\pm 1$ 
and we are in case \caseone. 
\item If $\lambda_2 \neq \pm 1$. Then $\lambda_2 \not\in \Q$, and by 
$(*)$, $Q(\lambda_2)=0$ (remark also that necessarily $deg_\Q\ \lambda_2=2$ and 
$deg_\Q\ \lambda_3=2$).

Denoting $\l=\lambda_2$, $det\ A=\pm 1$ gives 
$Spec_\C\ A = \{\l_1=\pm 1,\l_2=\l,\l_3=\pmlinv\}$. 
\begin{itemize}
\item If $\l\in\R$, we are in case \casetwoR. 
\item If $\l\in\C\setminus\R$, then $\l_3$ must equal $\bar{\l}_2$ (as 
complex roots of the real polynomial $P_A$). As $\l_3=\frac{-1}{\l}$ 
is impossible (it would imply $\l \bar{\l}=-1$), we get  $\l_3=\linv$, 
which implies $\l \bar{\l}=1$, so we are in case \casetwoC.
\end{itemize}
\end{itemize}
\item If $\deg\ P_{\l_1}=2$. Then $\lambda_1\not\in\Q$. As $P_A(\lambda_1)=0$, 
one can write $P_A=P_{\l_1}R$ with $R\in \Q[X]$ and $deg\ R=1$. Therefore 
$P_A$ has a root in \Q, which must be $\pm 1$. The same arguments as above 
show that we are in cases \casetwoC\ or \casetwoR, and that 
$deg_\Q\ \lambda_i \leq 2$ for $i=1\dots 3$.

\item If $\deg\ P_{\l_1}=3$. Then $deg_\Q\ \lambda_i = 3$ for $i=1\dots 3$ because 
none of the above cases can occur (there we had $deg_\Q\ \lambda_i \leq 2$ for $i=1\dots 3$). So we are in case \casethree. 
\end{itemize}
\end{proof}

{\bf Lemma \ref{few_multiple_roots} p.\pageref{few_multiple_roots}:}
\begin{proof}
Suppose that \l\ is a multiple root 
of $P_A$.

Assume $A$ falls in case \casethree. 
Then as $P_A'(\l)=0$ and $deg\ P_A'=2$, we get that $deg_\Q\ \l \leq 2$, 
which is impossible.

When $A$ falls in case \casetwoC\ or \casetwoR, it is easy to check 
that the condition $\l \neq \pm 1$ prevents $P_A$ from having 
a multiple root. Therefore $A$ falls in case \caseone. 
\end{proof}


{\bf Lemma \ref{proj_Stein_Reinhardt} p.\pageref{proj_Stein_Reinhardt}:}
\begin{proof}
This lemma is proved, as a particular case, in the next proof. 
\end{proof}

{\bf Lemma \ref{lowering_dimension} p.\pageref{lowering_dimension}:}
\begin{proof}
Up to dividing by some integer, we can assume that $v$ is unimodular. 
The matrix $M$ given by Lemma \ref{SLn_basis} yields an algebraic 
automorphism $h$ of $(\Cstar)^n$ by the usual correspondence. 
Now, up to conjugation by $h$, we can assume that any $g\in G_{struct}(E)$ 
is of the form:
 $$g(z_1, \dots, z_n)=(\alpha_1 z_1 z_2^*\cdots z_n^*,\  
\alpha_2 z_2^*\cdots z_n^*,\dots\dots,\ \alpha_n z_2^*\cdots z_n^*).$$ 
Therefore the $S^1$-action on $D$: 
$$\begin{array}[t]{c@{}c@{}c}
            S^1\times D                & \rightarrow & D \\
\Big( \theta\ ; (z_1, \dots, z_n) \Big) & \mapsto     & (\theta z_1, \dots, z_n)
         \end{array}
$$
commutes with the structural group $G_{struct}(E)$, 
so we get an induced $S^1$-action 
on $E$.

\smallskip
On the other hand, the projection map 
$$\applic{\pi}{D}{(\Cstar)^{n-1}}
{(z_1, \dots, z_n)} {(z_2, \dots, z_n)}$$ is 
$G_{struct}(E)$-invariant, with the action on $\pi(D)$ 
being {\em defined} by $g\Big(\pi(z_1, \dots, z_n)\Big) = 
\pi\Big( g(z_1, \dots, z_n) \Big)$. 

Therefore $\pi:D\rightarrow D'= \pi(D)$ induces a bundle map $q$:
$$\xymatrix{
\frac{\tilde{B}\times D}{\pi_1(B)}=& 
E \ar[rr]^{q} \ar[dr]   &   & E' \ar[dl]    & = \frac{\tilde{B}\times D'}{\pi_1(B)}\\
& & B &               
}$$

\smallskip
Now consider the following equivalence relation on $E$: 
$x \sim y$ when $f(x)=f(y)$ for all $S^1$-invariant  
holomorphic functions on $E$. Then by \cite{Heinz} we 
know that $E/\sim$ is a complex space, which is Stein 
if $E$ is Stein. 

We now prove that $q:E\rightarrow E'$ is a realization of that quotient, 
i.e., $x\sim y$ if and only if $q(x)=q(y)$. 
\begin{itemize}
\item If $q(x)\neq q(y)$ then as $E'$ is holomorphically 
separable (cf. section \ref{intro}), there exists $f\in \O(E')$ 
such that $f(q(x))\neq f(q(y))$. Then $f\circ q \in \O(E)$ 
separates $x$ and $y$. By construction, the $S^1$-orbits 
are contained in the fibers of $q$, so $f\circ q$ is $S^1$-invariant. 
Thus $x\not\sim y$.
\item If $q(x) = q(y)$ then $x$ and $y$ both belong to the same 
fiber of $q$, which is a one-dimensional annulus (connectedness 
of that fiber follows from the convexity of 
$log\ \big( D\cap (\Cstar)^n \big)$).\\ 
Any $S^1$-invariant  
function on $E$ is constant on each circle of that annulus, hence 
constant on the whole annulus (by the isolated zeros theorem). 
Thus $x\sim y$.
\end{itemize}

We still need to show that if $E'$ is Stein, then so is $E$. 
Assume $E'$ is Stein. Extending the $G_{struct}(E)$-action 
from $D$ to 
$\hat{D}= \{ (z_{1}, \ldots, z_{n})\ :\ 
z_{1}\in\Cstar,  (z_{2}, \ldots, z_{n})\in D'\}$, 
we get an equivariant map $D \hookrightarrow \hat{D}$, and (therefore) a 
corresponding injective bundle map: 

$$\xymatrix{
E\ \ar@{^{(}->}[rr] \ar[dr]   &   & \hat{E} 
                                     = \frac{\tilde{B}\times \hat{D}}{\pi_1(B)}\ar[dl] \\
                            & B &               
}$$
As $\hat{E}$ is a $\Cstar$-principal bundle over $E'$, it is Stein by 
the theorem of Matsushima and Morimoto. By local triviality of our bundles, 
$E$ is locally Stein in $\hat{E}$. Therefore $E$ is Stein by the Docquier-Grauert 
theorem. 
\end{proof}

{\bf Lemma \ref{psi_iso} p.\pageref{psi_iso}:}
\begin{proof}
Take $A_g$ and $A_{g'}$ in $\A(D)$. It is easy to check that 
$A_{gg'}$ is just $A_g A_{g'}$. The inverse matrix $A_g\inv$ 
comes from $g\inv$ by the above remark. Therefore $\A(D)$ is a group 
and $\psi$ is a morphism, which is onto by definition of $\A(D)$. 

Moreover, $\psi$ is one-to-one:\\
Assume $A_g = A_{g'}$. Then 
$\begin{array}[t]{r@{\hspace{1.5mm}}c@{\hspace{1.5mm}}l}
f_g\inv f_{g'} (p) & = & f_g\inv(A_{g'}p+b') \\
                   & = & A_g\inv A_{g'}p + A_g\inv b' -A_g b \\
                   & = & p + A_g\inv (b'-b).
\end{array}$\\
If $b'-b \neq 0$, then $A_g\inv (b'-b) \neq 0$, so 
$f_g\inv f_{g'}$ is a non-trivial translation of $\R^3$ 
which induces an automorphism of \Cvx. This is impossible 
because \Cvx\ contains no affine line. 

Thus $b'=b$. That proves the injectivity of $\psi$. 
\end{proof}

{\bf Lemma \ref{finite_index_subgroup} p.\pageref{finite_index_subgroup}:}
\begin{proof}
As $E$ is flat, it is given by a morphism 
$\rho :\pi_1(B)\rightarrow Aut(D)$ 
(see p.\pageref{E_given_by_representation}). 
Take $G:=\rho\inv(H)$. Then the bundle 
$$E'=(\tilde{B}\times D)/G \xrightarrow{D}B$$ 
has $G_{struct}(E')=H$ and is Stein if and only if $E$ is Stein, 
because it is a finite cover of $E$.
\end{proof}

{\bf Lemma \ref{if_spec_finite} p.\pageref{if_spec_finite}:}
\begin{proof}
We know by \cite{Shim} that any $g\in G_{struct}(E)$ has the form 
$$(*)\quad g(z_{1},\ldots, z_{n})=
(\alpha_{1}z_{1}^{a_{11}}\ldots z_{n}^{a_{1n}}, \ldots\ldots, 
\alpha_{n}z_{1}^{a_{n1}}\ldots z_{n}^{a_{nn}})$$ 
with $(a_{ij})_{1\leq i,j \leq n} \in \A(E)\subset GL_n(\Z)$ 
and $(\alpha_{1},\ldots,\alpha_{n})\in (\Cstar)^n$. 

As $Spec_\C\ \A(E)$ is finite, is follows from lemma 
\ref{Burnside}
that either $\A(E)$ is finite, or there exists an $\A(E)$-stable 
\Q-subspace $\{0\}\subsetneq V_1\subsetneq \Q^n$. In the latter case, 
call $\A_1\subset GL(V_1)$ the group of linear automorphisms induced 
by $\A(E)$. Being contained in $Spec_\C\ \A(E)$, $Spec_\C\ \A_1$ is 
finite, so we can apply Burnsides' result again to $\A_1$ acting on $V_1$. 
We repeat that procedure until we get a positive-dimensional subspace 
$V_k \subsetneq \Q^n$ such that $\A_k$ is finite (notice that if $dim\ V_k=1$, 
then finiteness of $Spec_\C\ \A_k$ implies finiteness of $\A_k$). 
Then a finite index subgroup of $\A(E)$ (namely, the kernel of 
$\A(E)\rightarrow \A_k$) fixes $V_k$ pointwise. 
By Lemma \ref{finite_index_subgroup}, we can assume that 
$\A(E)$ fixes $V_k$ pointwise.\\
Pick $v\in V_{k}\cap \Z^n$. Lemma \ref{lowering_dimension} 
applies and yields a bundle $E'\xrightarrow{D'}B$ with:\\ 
--- $E'$ is Stein if and only if $E$ is,\\
--- $dim\ D'=n-1$ and $D'$ is Stein by Lemma \ref{proj_Stein_Reinhardt},\\ 
--- any $g\in G_{struct}(E')$ still has the form (*) (with $n-1$ variables),\\
--- $Spec_\C\ \A(E')$ is a subset of $Spec_\C\ \A(E)$, and is therefore finite.

\smallskip
Even though $D'$ may be unbounded, the above properties allow 
us to repeat the procedure from the beginning. This sets up an induction on $n$, 
so we only need to prove the theorem for $n=1$. But in that case, as the fiber 
is of dimension one, $E$ is Stein by \cite{Mok}. 
\end{proof}

{\bf Lemma \ref{if_spec_subset_S1} p.\pageref{if_spec_subset_S1}:}
\begin{proof}
Assume $Spec_\C\ \A\subset S^1$. We prove finiteness of the set 
of all characteristic polynomials of all matrices in \A.

Take any $A\in\A$, and let $P$ be its characteristic polynomial. 
By assumption, all roots of $P$ have modulus $1$. As $deg\ P=n$, 
$P$ has $n+1$ coefficients. Each of them is an integer of 
modulus not greater than $n$, because it is given 
by an elementary symmetric expression in the roots of $P$. 
Therefore the set of polynomials has cardinal at most $(n+1)^{2n+1}$.
 
Conclusion:  $Spec_\C\ \A$ is finite. 
\end{proof}

\bigskip \noindent
{\small {\bf Acknowledgements.} The second named author wishes to 
warmly thank A. Huckleberry for his hospitality in Bochum in 2001, where part 
of this work was done. We thank as well S. Chase, I. Chatterji, J. Hubbard, P. Liardet 
and K. Wortman for useful conversations. 
This work has been done independently from \cite{Pfl-Zwo}.
}


\bigskip
        \begin{thebibliography}{2}
            \bibitem [C\oe-L\oe]{Coe-Loe} 
                       {{\bf G. C\oe ur\'e, J.-J. L\oe b:}
                        ``A counterexample to the Serre problem 
                        with a bounded domain of $\C^2$ as fiber''. 
                        Ann. of Math., 122 (1985), 329--334.}
            \bibitem [Die-For]{Die-For}
                       {{\bf K. Diederich, J. Fornaess:}
                        ``Pseudoconvex domains: bounded strictly 
                        plurisubharmonic exhaustion functions''.
                        Inventiones Math. 39, 129--141 (1977).}
            \bibitem [Heinz]{Heinz}
                       {{\bf P. Heinzner:}
                        ``Kompakte Transformationsgruppen Steinscher 
                        R\"aume'' [Compact transformation groups 
                        of Stein spaces]. 
                        Math. Ann. 285 (1989), no. 1, 13--28.} 
	    \bibitem[Hell]{Hell}{\bf Y. Hellegouarch:}
	                Invitation to the mathematics of Fermat-Wiles. 
                        Academic Press, Inc., San Diego, CA, 2002.
            \bibitem [Hir]{Hirsch}
                       {{\bf A. Hirschowitz:}
                        ``Domaines de Stein et fonctions holomorphes born\'ees''. 
                        Math. Ann.  213  (1975), 185--193.} 
            \bibitem [Kaup]{Kaup}
                       {{\bf W. Kaup:}
                        ``Reelle Transformationsgruppen und invariante Metriken 
                        auf komplexen R\"aume''. 
                        Invent. Math.  3  (1967), 43--70.} 
            \bibitem [Mok]{Mok}
                       {{\bf N. Mok:}
                        ``Le probl\`eme de Serre pour les surfaces de Riemann''. 
                        C. R. Acad. Sci. Paris S\'er. A-B  290  (1980), 
                        no. 4, A179--A180.} 
            \bibitem [Nara]{Nara}
                       {{\bf R. Narasimhan:}
                        ``Several complex variables''. 
                        Chicago Lectures in Mathematics. 
                        University of Chicago Press, Chicago, IL, 1995.} 
            \bibitem [Pfl-Zwo]{Pfl-Zwo}
                       {{\bf P. Pflug, W. Zwonek:} 
                       ``The Serre problem with Reinhardt fibers''. 
                       Ann. Inst. Fourier, 54, 1 (2004), 129-146.}
            \bibitem [Ser]{Ser}
                       {{\bf J.-P. Serre:}
                        ``Quelques probl\`emes globaux relatifs aux
                        vari\'et\'es de Stein''. 
                        Colloque sur les fonctions de plusieurs 
                        variables, Bruxelles, 1953, 57--68.}
            \bibitem [Shim]{Shim}
                       {{\bf S. Shimizu:}
                        ``Automorphisms of bounded Reinhardt domains''. 
                        Japan. J. Math. (N.S.)  15  (1989),  no. 2, 385--414.}
            \bibitem [Sib]{Sib}
                       {{\bf N. Sibony:}
                        ``Fibr\'es holomorphes et m\'etrique de Carath\'eodory''.
                        C. R. Acad. Sci. Paris S\'er. A  279  (1974), 261--264.} 
            \bibitem [Siu 1]{Siu1}
                       {{\bf Y.T. Siu:}
                        ``All plane domains are Banach-Stein''. 
                        Manuscripta Math.  14  (1974), 101--105.} 
            \bibitem [Siu 2]{Siu2}
                       {{\bf Y.T. Siu:}
                        ``Holomorphic fiber bundles whose fibers are 
                        bounded Stein domains with zero first Betti number''. 
                        Math. Ann. 219 (1976), no. 2, 171--192.}
            \bibitem [Siu 3]{Siu3}
                       {{\bf Y.T. Siu:}
                        ``Pseudoconvexity and the problem of Levi''.
                        Bull. Am. Math. Soc. 84 (1978), 481--512.}
            \bibitem [Sko]{Skoda}
                       {{\bf H. Skoda:}
                     ``Fibr\'es holomorphes \`a base et \`a fibre de Stein''. 
                        Invent. Math. 43 (1977), no. 2, 97--107.}
            \bibitem [Steh]{Steh}
                       {{\bf J.-L. Stehl\'e:}
                      ``Fonctions plurisousharmoniques et convexit\'e 
                      holomorphe de certains fibr\'es analytiques''. 
                      C. R. Acad. Sci. Paris S\'er. A 279 (1974), 235--238.}
           \bibitem [Stein]{Stein}
                       {{\bf K. Stein:}
                        ``\"Uberlagerungen holomorph-vollst\"andiger komplexer 
                         R\"aume''.  
                        Arch. Math.  7  (1956), 354--361.}
            \bibitem [Steinb]{Steinb}
                       {{\bf R. Steinberg:}
                        ``Conjugacy classes in algebraic groups''.
                        Notes by Vinay V. Deodhar. 
                        Lecture Notes in Mathematics, Vol. 366.
                        Springer-Verlag, Berlin-New York, 1974.} 
            \bibitem [Zaf]{Zaf}
                       {{\bf D. Zaffran:}
                        ``Serre problem and Inoue-Hirzebruch surfaces''.
                        Math. Ann.  319  (2001),  no. 2, 395--420.} 

            \bibitem [Zwo]{Zwo}
                       {{\bf W. Zwonek:}
                        ``On hyperbolicity of pseudoconvex Reinhardt domains''.
                        Arch. Math. (Basel)  72  (1999),  no. 4, 304--314.} 
\end{thebibliography}
\end{document}